\documentclass[10pt]{smfart}
\usepackage{a4wide}
\usepackage{amsmath,amssymb,smfthm}
\usepackage[francais]{babel}
\usepackage[latin1]{inputenc}
\usepackage[all]{xy}

\def\L{\mathcal{L}}
\def\La{\Lambda} 
\def\la{\lambda}

\def\O{\mathcal{O}}
\def\F{\mathcal{F}}

\def\*{^\times }
\def\dpt{\displaystyle}
\def\G{\mathcal{G}}

\def\a{\alpha}

\def\ph{\varphi}

\def\limpl{\Longrightarrow}
\def\drt{\rightarrow}
\def\ldrt{\longrightarrow}

\def\Ql{\mathbb{Q}_\ell}
\def\Qlb{\overline{\mathbb{Q}}_\ell}

\def\Z{\mathbb{Z}}
\def\N{\mathbb{N}}
\def\Zl{\mathbb{Z}_\ell}
\def\Hom{\text{Hom}}

\def\Rpsi{\text{R}\Psi_{\hspace{-1pt}\eta}}
\def\rpsi{\text{R}\Psi}
\def\Gal{\text{Gal}}
\def\={\! = \!}

\def\spec{\text{Spec}}

\def\limp{\underset{\longleftarrow}{\text{ lim }}\;}
\def\limi{\underset{\longrightarrow}{\text{ lim }}\;}

\def\iso{\xrightarrow{\;\sim\;}}

\def\End{\text{End}}

\def\xrig{\xrightarrow}

\def\Ext{\text{Ext}}

\def\AA{\mathcal{A}}

\def\<{<\hspace{-1mm}}
\def\>{\hspace{-1mm}>}

\def\dem{{\it D{\'e}monstration. }}

\def\Fil{\text{Fil}}

\def\Perv{\text{Perv}}
\def\et{\text{ét}}

\author{  Laurent Fargues} 
\address{CNRS-université Paris-Sud-IHES} 
\email{laurent.fargues@math.u-psud.fr} 
\date{}

\begin{document}

\title{Filtration de monodromie et cycles évanescents formels}
\maketitle

\begin{abstract}
D'après Berkovich \cite{Berk3}, Fujiwara \cite{Fuj2} et Huber \cite{Hu3} la fibre des cycles évanescents en un point de la fibre spéciale ne dépend que du complété formel en ce point. Nous raffinons ce résultat en démontrant l'invariance par complétion formelle de la filtration de monodromie perverse sur la fibre des cycles évanescents. Ce résultat est utilisé de façon cruciale par Boyer dans \cite{Boyer2}. 
\end{abstract}

\begin{altabstract}
Berkovich \cite{Berk3}, Fujiwara \cite{Fuj2} and Huber \cite{Hu3} 
have proved that the fiber of the vanishing cycles at a point of the special fiber depends only on the formal completion at this point. We refine this result 
and prove the invariance under formal completion of  the perverse monodromy filtration on the fiber of vanishing cycles. This result is used in an essential way by Boyer in \cite{Boyer2}.
\end{altabstract}

\section{Notations et énoncé du théorème}

Soit $S= \spec ( \O)$ un trait hensélien excellent. On note $k$
le corps résiduel de $\O$ que l'on suppose de caractéristique $p>0$
et $\pi$ une uniformisante de $\O$. On note $\eta$ le point générique
de $S$, $s$ son point fermé, $\overline{\eta}$ un point géométrique
au dessus de $\eta$, $\overline{S}=\spec (\overline{\O})$ la fermeture
intégrale de $S$ dans $\overline{\eta}$ et $\overline{s}$ le point
géométrique au dessus de $S$ associé. 

Soit $\ell$ un nombre  premier différent de $p$.
Pour $X$ un $S$-schéma de type fini et $\F\in \mathbb{D}^b_c (X,\Qlb)$
on note 
$$
\Rpsi \F = \overline{i}^* R \overline{j}_* \F 
\in \mathbb{D}^b_c (X_{\overline{s}}, \Qlb )
$$
le complexe des cycles proches de Deligne.

Si l'on muni $X_\eta$ et $ X_{\overline{s}}$ de la perversité intermédiaire
 et si $\F \in \Perv (X_\eta )$ est un faisceau pervers sur la fibre générique 
alors 
$$
\Rpsi \F\in \Perv (X_{\overline{s}})
$$
et est muni d'une action de $\Gal (\overline{\eta}|\eta)$. 
Si $I$ désigne le sous-groupe d'inertie de $\Gal (
\overline{\eta}|\eta)$ et $P$ son sous-groupe de ramification
sauvage 
$$
P/I \simeq  \prod_{ \ell'\neq p} \mathbb{Z}_{\ell'} (1)
$$
et si $T$ désigne un élément de $I$ s'envoyant sur un générateur
de $\Zl (1)$ dans l'isomorphisme précédent, 
$$
T\in \End_{\Perv (X_{\overline{s}})} ( \Rpsi \F)
$$
qui est quasi-unipotent au sens où si $n$ est suffisamment divisible
$T^n-Id$ est un endomorphisme nilpotent de $\Rpsi \F$. Cela permet
de définir l'endomorphisme nilpotent 
$$
N= \log T = \frac{1}{n} \log \left ( Id + ( T^n -Id) \right )
\text{  pour } n>>0
$$
$$
N : \Rpsi \F (1) \ldrt \Rpsi \F
$$
Il définit une filtration de monodromie dans $\Perv
(X_{\overline{s}})$, $\Fil^\bullet \Rpsi  \F$
qui ne dépends pas du choix de $T$, 
 et donc (cf. 
la section \ref{catfiltperv}) un objet $K (\F)$ de la catégorie 
dérivée filtrée $\mathbb{D}^b_c F (X_{\overline{s}},\Qlb)$.
 \\

Soit $x : \overline{s} \ldrt X_{\overline{s}}$ un point géométrique
fermé. Le complexe de $\Qlb$-espaces vectoriels $\Rpsi (\F)_x$ 
est muni d'une action de $N$ mais n'est pas pervers en général (en quelque sens
que ce soit). Néanmoins il provient par oubli de la filtration  
du complexe filtré $K (\F)_x \in \mathbb{D}^b F(\Qlb)$ la catégorie dérivée
filtrée des $\Qlb$-espaces vectoriels.

Rappelons 

\begin{theo}[Berkovich \cite{Berk3} th. 3.1, Fujiwara \cite{Fuj2} cor.
  7.1.7, Huber
  \cite{Hu3} prop. 3.15]
Soit $(X',\F',x')$ un triplet analogue à $(X,\F,x)$. Tout 
isomorphisme 
$$
f: \widehat{\O}_{X,x} \iso \widehat{\O}_{X',x'}
$$
entre les complétés formels des anneaux locaux
 tel que 
$$
\F_{ |\spec ( \widehat{\O}_{X,x}[\frac{1}{\pi}])} =
f^* \left ( \F'_{|\spec ( \widehat{\O}_{X',x'}[\frac{1}{\pi}])} \right )
$$
induit naturellement un isomorphisme 
$$
\Rpsi (\F)_x \iso \Rpsi (\F)_{x'}
$$

\end{theo}

Le but de cet appendice est d'établir la version raffinée suivante de
cet énoncé :

\begin{theo}
Soit $(X',\F',x')$ un triplet analogue à $(X,\F,x)$.
Soient $K (\F)$, resp. $ K(\F')$, le complexe filtré associé 
à $\rpsi (\F)$, resp. $\rpsi (\F')$, 
 muni de sa filtration de monodromie.
 
Tout 
isomorphisme 
$$
f: \widehat{\O}_{X,x} \iso \widehat{\O}_{X',x'}
$$
 tel que 
$$
\F_{ |\spec ( \widehat{\O}_{X,x}[\frac{1}{\pi}])} =
f^* \left ( \F'_{|\spec ( \widehat{\O}_{X',x'}[\frac{1}{\pi}])} \right )
$$
induit naturellement un isomorphisme 
$$
K (\F)_x \iso K (\F)_{x'}
$$
dans la catégorie dérivée filtrée des $\Qlb$-espaces vectoriels 
$\mathbb{D}^b F (\Qlb)$.
\end{theo}

La démonstration de ce théorème repose sur le théorème de
désingularisation de Popescu (cf. par exemple \cite{Pop}) ainsi que la
définition d'une catégorie de faisceaux pervers $\ell$-adiques sur des
schémas pour lesquels les théorèmes de finitude de Deligne de SGA
$4^{1/2}$ et l'existence d'un module dualisant 
 ne sont pas connus comme le complété formel d'un anneau
local d'un schéma de type fini sur un corps.
La nécessité des coefficients $\ell$-adiques est imposée par la
définition de l'opérateur de monodromie comme logarithme d'un autre
opérateur. Sans théorèmes de finitude l'argument d'autorité ``le cas
des coefficients de torsion s'étend aussitôt au cas $\ell$-adique'' ne
peut être employé ce qui induit de difficiles contorsions basées sur
les travaux d'Ekedahl (\cite{Eke}). 
\\
Le résultat de cet article est utilisé de façon cruciale par P.Boyer dans \cite{Boyer2}. 

L'auteur ce cet appendice tient à remercier Ofer Gabber pour
d'intéressantes discussions sur le sujet. 

\section{Catégories dérivées filtrées}

Soit $\mathcal{A}$ une catégorie abélienne. On utilise 
les catégories dérivées filtrées 
pour des filtrations décroissantes finies (et pas seulement
finies degré par degré) 
$\mathbb{D} F ( \mathcal{A}),\mathbb{D}^b F (\mathcal{A}), 
\mathbb{D}^+F (\mathcal{A}) $ 
telles qu'elles sont
définies dans le chapitre V de \cite{Illusie1}, cf. également le chapitre 3 de 
\cite{BBD}.
Elles sont munies de foncteurs d'oubli de la filtration, $Gr^i$ pour 
$i\in \Z$, et
$ \Fil^i/\Fil^j$ pour $i<j$ vers les catégories dérivées usuelles.
Les foncteurs $\Fil^i/\Fil^j$ se factorisent par les catégories
dérivées filtrées et on notera également par abus de notations
$\Fil^i/\Fil^j$ ces factorisations.

\section{Catégories dérivées filtrées et faisceaux pervers filtrés}
\label{catfiltperv}

On reprend les hypothèses du début de la section 3.1.5 de 
\cite{BBD}. Plus précisément, 
supposons que $\mathcal{A}$ possède suffisamment d'injectifs, 
que $\mathcal{D}$ soit une sous-catégorie triangulée pleine de $\mathbb{D}^+
(\mathcal{A})$ et que $( \mathcal{D}^{\leq 0}, \mathcal{D}^{\geq 0})$
soit une $t$-structure sur $\mathcal{D}$ de coeur $\mathcal{C}
= \mathcal{D}^{\leq 0} \cap  \mathcal{D}^{\geq 0}$.
Soit $\mathcal{D}F$ la sous-catégorie triangulée de $\mathbb{D}^+ F (\mathcal{A})$
formée des objets  $A$ tels que $\forall i\in \Z\;
Gr^i A \in \mathcal{D}$.

\begin{prop}\label{PervFilt}
Il y a une équivalence de catégories entre la catégorie formée des couples
$(X,  Fil^\bullet X)$ où $X\in \mathcal{C}$ et $Fil^\bullet X$
est une filtration décroissante finie de $X$  et la sous-catégorie
de $ \mathcal{D} F$ formée des objets
 $A$ tels que $\forall i \; Gr^i A \in \mathcal{C}$. 
Dans cette équivalence les foncteurs de graduation $Gr^i$ se correspondent et plus généralement
si $A$ est associé à $X$ 
$\forall i<j\;\; \Fil^{\,i}/\Fil^{\,j} (A) = Fil^i X/\Fil^{\,j} X$.  
Ainsi via le foncteur d'oubli de la filtration $\omega = \Fil^{-
\infty} /
\Fil^{\,\infty}
: \mathcal{D}F
\drt \mathcal{D}$, $\omega A = X$.
\end{prop}
\dem
Soit $A\in \mathcal{D}F$ comme dans l'énoncé. Il y a des triangles
exactes dans $\mathcal{D}$ pour $j>i$
$$
\xymatrix@R=4mm@C=4mm{
\Fil^i/\Fil^j (A) \ar[rr] && \Fil^{i-1}/\Fil^j (A) \ar[ld] \\
& Gr^{i-1} A \ar[lu]^{+1}
}
$$
desquels on déduit aisément en procédant 
à $j$ fixé par récurrence sur $i, i>j$ 
à partir de $i=j-1$, que
$$
\forall j>i\;\; \Fil^i /\Fil^j ( A) \in \mathcal{C}
$$ 
et en particulier $A\in \mathcal{C}$. Donc, dans le triangle exacte 
$$
\xymatrix@R=4mm@C=4mm{
\Fil^i A \ar[rr] && \omega A \ar[ld] \\
& \Fil^{-\infty}/\Fil^i A \ar[lu]^{+1}
}
$$
les trois objets sont dans $\mathcal{C}$. On en déduit que 
$\Fil^i A \drt \omega A $ est un monomorphisme dans $\mathcal{C}$
et définit une filtration finie  de $X=\omega A$ ayant comme gradués
les $Gr^\bullet A$. Cela définit un foncteur $A\mapsto (X,
\Fil^\bullet)$.
Montrons qu'il est pleinement fidèle.

Commençons par remarquer que pour $A,B$ dans $\mathcal{D}F$ comme précédemment 
\begin{eqnarray}\label{anolt}
\Hom^{-1}_{\mathcal{D}F} (A,B)=0
\end{eqnarray}
Utilisons pour cela la suite spectrale (3.1.3.5) p.78 de \cite{BBD} :
$$
E^{pq}_1 = \left \{ \dpt{\prod_{j-i=p} \Hom^{p+q} ( Gr^i A, Gr^j B)
\;\text{ si } p \geq 0} \atop
 0\; \text{ si } p<0 \right.
$$
$$
E_1^{pq} \limpl \Hom_{\mathcal{D}F}^{p+q} (A,B)
$$
Les $ Gr^i A$ et $Gr^j B$ étant dans $\mathcal{C}$, $ E^{pq}_1 =0$
 si $p+q<0$ et l'égalité (\ref{anolt}) en résulte.
\\
Notons $(X,\Fil^\bullet X), (Y,\Fil^\bullet Y)$ les objets associés à $A$ et $B$.
On pourrait être tenté d'utiliser la suite spectrale précédente 
couplée à une même suite spectrale pour $\mathbb{D}^+F (\mathcal{C})$  
pour montrer que 
$\Hom (A,B)\iso \Hom ((X,\Fil^\bullet X), (Y,\Fil^\bullet
Y))$. Néanmoins cela s'avère délicat à cause de l'identification de
certaines flèches dans la suite spectrale et parce qu'à priori une telle suite
spectrale n'existe pour les objets filtrés de $\mathcal{C}$ que si
$\mathcal{C}$  possède suffisamment d'injectifs (néanmoins ce type
d'argument est sûrement possible). 
Supposant $B\neq 0$, procédons plutôt par récurrence sur l'amplitude de
la filtration $\text{sup} \{ i\;|\; Gr^i B\neq 0 \}- \text{inf} \{
i\;|\; Gr^i B\neq 0 \}$. Lorsque cet entier est égal à $0$, soit $i$
tel que $Gr^i B\neq 0$. On a donc $\Fil^{i+1} B =0$ ainsi que $\Fil^i B =B$ ce
qui implique $\Fil^{i+1} X=0$ et $\Fil^i X=X$. Alors
$$
\Hom_{\mathcal{D} F} (A,B) = \Hom_{\mathcal{D}} (\Fil^{-\infty}/\Fil^{\, i+1} A, B) = \Hom (X/\Fil^{i+1} X, Y) = \Hom ((X,\Fil^\bullet X), (Y,\Fil^\bullet Y))
$$
Supposons maintenant l'amplitude de la filtration non-nulle et soit
$i$ tel que $\Fil^{i} B\neq 0$, $\Fil^{i +1} B=0$ et donc $\Fil^{i+1}
Y =0$.  
D'après l'annulation (\ref{anolt}) il y a un morphisme de suite exactes 
$$
\xymatrix@C=4mm{
0 \ar[r] & \Hom_{\mathcal{D} F} (A, \Fil^i B) \ar[r]\ar[d]^u & \Hom_{\mathcal{D} F} (A,B) \ar[r]\ar[d]^v & 
\Hom_{\mathcal{D} F} (A, \Fil^{-\infty}/\Fil^i B) \ar[r]\ar[d]^w & \Hom^1_{\mathcal{D}F} (A, \Fil^i B)\ar[d]^x \\
0 \ar[r] & \Hom (\Fil^\bullet X, \Fil^i Y ) \ar[r] & 
\Hom (\Fil^\bullet X, \Fil^\bullet Y) \ar[r] & \Hom (\Fil^\bullet X, \Fil^\bullet Y/\Fil^i Y) \ar[r]^(.55)\delta & \Ext^1 (\Fil^\bullet X,  \Fil^i Y) 
}
$$
où les $\Hom$ de la ligne du bas sont des morphismes d'objets
filtrés, $\Fil^i Y$ étant muni de la filtration induite, et 
où  $\Ext^1 (\Fil^\bullet X,  \Fil^i Y) $ est défini comme
$$\Ext^1_{\mathcal{C}, \text{Yoneda}} ( X/\Fil^{i+1} X, \Fil^i Y)$$
De plus si $[\xi]$ désigne la classe de l'extension 
$$
\xi \;: \;\;\;\;\; 0 \drt \Fil^i Y \drt Y \drt Y/\Fil^i Y\drt 0
$$
alors l'application $\delta$ est définie par le composé 
$$
\xymatrix@R=5mm@C=5mm{
 \Hom (\Fil^\bullet X, \Fil^\bullet Y/\Fil^i Y) \ar[r]^(.44)\sim &  \Hom
 (\Fil^\bullet X/\Fil^{i+1} X , \Fil^\bullet Y/\Fil^i Y) \ar[d] \\
&  \Hom
 ( X/\Fil^{i+1} X, Y/\Fil^i Y) \ar[r] & \Ext^1_{\mathcal{C}, \text{Yoneda}} (
 X/\Fil^{i+1} X, Y) \\
& f \ar@{|->}[r] & f^*[\xi]
}
$$
Avec cette définition on vérifie facilement que le carré de droite est
commutatif dans le diagramme précédent (pour la définition de
l'application $x$ cf. la ligne qui suit). 
 \\ 
D'après le cas étudié précédemment initialisant la récurrence l'application $u$ est un isomorphisme. Par hypothèse de récurrence $w$ en est également un. Quant à l'application $x$ cela résulte de 
$$
\Hom^1_{\mathcal{D}F} (A, \Fil^i B) = \Hom^1_{\mathcal{D}} (
\Fil^{-\infty}/\Fil^{i+1} A,  \Fil^i B) = \Ext^1_{\mathcal{C}, \text{Yoneda}} ( X/\Fil^{i+1} X, Y)
$$
L'application $v$ est donc un isomorphisme.

La surjectivité essentielle du foncteur $A\mapsto (X,\Fil^\bullet X)$ résulte de ce que tout morphisme de
complexes dans $\mathbb{D}^+ (\mathcal{A})$ peut être représenté 
par un morphisme injectif de complexe d'objets injectifs
(toute flèche dans $\mathbb{D}^+ (\mathcal{A})$ peut être vue comme
un cône d'une autre flèche) et que donc si $(X, \Fil^\bullet X)$
est comme dans l'énoncé, on peut représenter $\Fil^i X$par un complexe 
$(K^{ij})_{j\in \Z}$ et les flèches $\Fil^{i+1} X \drt \Fil^i X
$ par des morphismes injectifs $K^{i+1\bullet}\hookrightarrow 
K^{i\bullet}$. L'objet de $\mathcal{D}F$ associé convient alors.
\qed

\section{Faisceaux pervers sans conditions de finitude d'après
Gabber}\label{tGabber}

Le théorème suivant est un cas particulier des résultats 
de \cite{Gabber1}. On y note $Ri_x^! \F^\bullet $ la fibre en un point
géométrique au dessus du point générique de $\overline{\{ x \} }$ de
$R i_{\overline{\{ x \} }}^! \F^\bullet$. Ainsi, si $\overline{x}\drt
X$ est un tel point  géométrique 
$$
\def\objectstyle{\scriptstyle}
\def\labelstyle{\scriptstyle}
\mathcal{H}^i (i_x^! \F^\bullet ) = \underset{\xymatrix@-1.4pc{ &
    U\ar[d]^{a\text{\tiny  ,étale}} \\
    \overline{x}\ar[r]\ar[ru] & X }}{\limi} H^i_{a^{-1}\left (\overline{\{
    x \}}\right )} (U,\F^\bullet)
$$
L'annulation d'un tel groupe ne dépend pas du choix de $\overline{x}$
au dessus de $x$, de même pour $\mathcal{H}^i (i_x^*\F^\bullet)$.

\begin{theo}[Gabber]
Soit $X$ un schéma noethérien de dimension finie. Soit $\La$ 
un anneau. Posons 
$$
\,^p\mathbb{D}^{\leq 0} (X,\La) = \{ \F^\bullet
\in\mathbb{D} (X_{\text{ét}} ,\La) \; | \;
\forall x\in X\; \forall i>-\dim \overline{\{ x\}} \;\;
\mathcal{H}^i ( i_x^* \F^\bullet ) = 0 \; \}
$$
$$
\,^p\mathbb{D}^{\geq 0} (X,\La) = \{ \F^\bullet \in\mathbb{D}^+ (X_{\text{ét}}
,\La)  \; | \;
\forall x\in X\; \forall i<-\dim \overline{\{ x\}} \;\;
\mathcal{H}^i ( R i_x^! \F^\bullet ) = 0 \; \}
$$
Alors, $(\,^p\mathbb{D}^{\leq 0},\,^p\mathbb{D}^{\geq 0})$ définit une
t-structure sur $\mathbb{D} (X_{\text{ét}}
,\La)$ et induit une t-structure sur $\mathbb{D}^+ (X_{\text{ét}}
,\La)$, $\mathbb{D}^- (X_{\text{ét}}
,\La)$ et $\mathbb{D}^b (X_{\text{ét}}
,\La)$.

On note $\text{Perv} (X,\La)= \,^p\mathbb{D}^{\leq 0} \cap
\,^p\mathbb{D}^{\geq 0} \subset \mathbb{D}^b (X_{\text{ét}},\La)$, une
catégorie abélienne $\La$-linéaire.

Si de plus $\La=\Z/n\Z$, $n$ est inversible sur $X$ et $X$ possède un
module dualisant au sens de Grothendieck (SGA 5, exposé I) alors 
la t-structure précédente induit une t-structure sur $\mathbb{D}^b_c (X_{\et},\La)$
et tout objet de 
$\text{Perv} (X) \cap  \mathbb{D}^b_c (X_{\text{ét}},\La)$ est de longueur finie.
\end{theo}

L'intérêt de ce théorème est qu'il permet de construire la
t-structure perverse sans hypothèse de finitude  sur les $R
j_*$ pour $j$ une immersion ouverte, $R i^!$ pour $i$ une immersion
fermée, ou bien l'existence d'un module dualisant
 (ce qui est le cas considéré dans
\cite{BBD} que l'on retrouve  à la fin de l'énoncé du théorème
précédent). La
dernière hypothèse sur l'existence d'un module dualisant permet de s'assurer que les opérateurs de
troncature perverse préservent les faisceaux pervers à cohomologie 
bornée constructible.

\begin{exem}
 Supposons $X$ excellent régulier de dimension $d$
et $\L$ un faisceau localement constant sur $X_{\et}$. Alors, $\L[d]\in \text{Perv} (X,\La)$. C'est une conséquence de ce que $\forall x\in X\;$ le schéma réduit  sous-jacent à $\overline{\{ x\}}$ est génériquement régulier et du théorème principal de \cite{FujGabber}. 
\end{exem}

\section{Extension aux coefficients $\ell$-adiques}

Soit $L_\la |\Ql$ une extension de degré fini, $\O_\la$ son
anneau des entiers et $\varpi_\la$ une uniformisante de $\O_\la$.

Soit $X$ un schéma. Nous allons utiliser le formalisme de
\cite{Eke}. Certains des énoncés de \cite{Eke} sont trop elliptiques
(par exemple le (v) du théorème 3.6 où la t-structure tautologique
n'est pas définie). De plus nous aurons besoin de généralisations,
c'est pourquoi nous rappelons le formalisme de \cite{Eke} dans les
paragraphes qui suivent.

\subsection{Généralités : $R\pi_*$ et $ \mathbb{L}\pi^*$}

 Soit
$\widetilde{X}_\et^\N-\O_{\la\bullet}$ la catégorie abélienne des systèmes
projectifs $(\F_n)_{n\geq 1}$
$$
\drt \F_{n+1} \drt \F_n \drt \dots \drt \F_1
$$
où $\F_n$ est un $\O_\la/\varpi_\la^n$-module sur $X_\et$ 
le site étale de $X$ et les
flèches de transition sont $\O_\la$-linéaires.  Soit
$\widetilde{X}_\et-\O_\la$ la catégorie des $\O_\la$-modules sur
$X_\et$. Soit 
$$
 (\pi^*,\pi_*) : \widetilde{X}_\et^\N-\O_{\la\bullet} \ldrt \widetilde{X}_\et-\O_\la
$$
le couple de foncteurs défini par 
$$
\pi^* \F = \left ( \F/\varpi_\la^n \F \right )_{n\geq 1}
$$
et 
$$
\pi_* (\F_n)_n = \underset{n}{\limp} \F_n
$$
On note donc $\pi^*$ pour le foncteur image inverse associé au morphisme de topos annelés et $\pi^{-1}$ le
foncteur ``sans anneau'' :
$$
\pi^* \F= \pi^{-1}\F\otimes_{\pi^{-1} (\O_\la)}\O_{\la\bullet}
$$ 
\begin{rema}
On remarquera qu'en général $R\pi_*$ n'est pas de dimension
cohomologique finie (en général on a seulement $\text{cd} (R\pi_*) \leq
\text{cd}_\ell (X_{\text{ét}}) + 1$).
Néanmoins $\mathbb{L}\pi^*$ est de dimension cohomologique finie
(égale à $1$).
\end{rema}

\begin{defi}
On note  $\mathbb{D}(X^\N,\O_{\la \bullet})$  pour la catégorie
dérivée $ 
\mathbb{D} (\widetilde{X}_\et^\N-\O_{\la\bullet})$. 
\end{defi}

\begin{rema}
Le morphisme de topos annelés $\widetilde{X}_{\et}^\N \ldrt
\widetilde{X}_{\et}$ est un cas particulier de morphisme associé à un
topos fibré tel qu'expliqué dans \cite{Illusie3}. Nous renvoyons en
particulier à \cite{Illusie3} pour les généralités concernant les
résolutions dans les topos fibrés. 
\end{rema}

\subsection{Complexes essentiellement nuls/constants}

\begin{defi}[\cite{Eke} déf. 2.1]
Un système projectif $(\F_n)_n$ est dit essentiellement nul si $\forall n$,
localement sur $X_\et$, $\exists N\;\; \F_{n+N} \drt \F_n$ est le
morphisme nul.

Un complexe $C\in \mathbb{D}(X^\N,\O_{\la \bullet})$
est dit essentiellement nul si $\forall i\; \mathcal{H}^i (C)$ l'est.
\end{defi}

Remarquons que lorsque $X$ est quasicompact on peut
enlever ``localement sur $X_\et$'' dans la définition précédente.

\begin{rema}\label{reth}
Les systèmes essentiellement nuls forment une sous-catégorie
épaisse de la catégorie des systèmes projectifs de faisceaux 
et la catégorie quotient des systèmes projectifs de faisceaux par les systèmes
essentiellement nuls est donc abélienne. 
 On peut donc définir le localisé de $\mathbb{D}
(X^\N,\O_{\la \bullet})$ par rapport au flèches de cône
essentiellement nul. C'est une catégorie triangulée. En effet, si
$\mathcal{D}$ est une catégorie triangulée et $H$ est un foncteur
cohomologique sur $\mathcal{D}$ (à valeurs dans une catégorie
abélienne, par définition) alors les flèches $f$ dont un cône $C(f)$
vérifie $\forall i \; H( C(f)[i])=0$ permettent un calcul des fraction
à gauche et à droite, et la catégorie localisée est triangulée
(appliquer cela à $\mathcal{D} = \mathbb{D} (X^\N,\O_{\la \bullet})$ et $H$
le foncteur composé 
$\mathcal{D} \xrig{\; \mathcal{H}^0\;} \widetilde{X}_\et^\N 
\ldrt \widetilde{X}_\et^\N/\text{ess. nuls}$). 
\end{rema}

\begin{defi}[\cite{Eke} section 1]
On note $\mathbb{D}^{e+} (X^\N,\O_{\la \bullet})$ la sous-catégorie de
$\mathbb{D}(X^\N,\O_{\la \bullet})$ formée des complexes $C$ tels que
pour $i<<0\;\; \mathcal{H}^i (C)$ soit essentiellement nul. Il s'agit
donc des complexes isomorphes aux objets de $\mathbb{D}^+$
dans la catégorie modulo les complexes
essentiellement nuls de la remarque précédente.
\end{defi}

Si $C\in\mathbb{D}^+ (X^\N,\O_{\la\bullet})$ est essentiellement nul alors $R\pi_* C=0$.
On en déduit que si $C\in \mathbb{D}^{e+} (
X^\N,\O_{\la\bullet})$ le système projectif
$$
\tau (C) = (\dots \drt \tau_{\geq i} C\drt \tau_{\geq i+1} C\drt \dots )
$$
est tel que le système projectif $R\pi_*\tau (C):=
(R\pi_* \tau_{\geq i} C)_i$ est essentiellement constant. Cela permet
d'étendre $R\pi_*$ à $ \mathbb{D}^{e+} (
X^\N,\O_{\la\bullet})$ en posant 
$$
R\pi_* C = \limp R\pi_* \tau (C)
$$

\begin{defi}
Un objet de $\mathbb{D}^+ (X^\N,{\O_{\la\bullet}})$ est
essentiellement constant s'il est isomorphe modulo les complexes
essentiellement nuls à un objets de la forme $\pi^{-1} (D)$.
\end{defi}

\begin{rema}
Le complexe $C$ est essentiellement constant ssi l'application
d'adjonction $\pi^{-1} R\pi_* C\drt C$ a un cône essentiellement
nul ssi $\forall i$ l'application $\pi^{-1}\pi_* \mathcal{H}^i (C) \drt
\mathcal{H}^i (C)$ a un noyau et conoyau essentiellement nuls.
\end{rema}

\subsection{Le foncteur $\pi^* (
  \O_\la/\varpi_\la) \overset{\mathbb{L}}{\otimes}_{\O_{\la\bullet}} -$}

 Bien que $\pi^* (\O_\la/\varpi_\la)$ ne
soit pas un $\O_{\la \bullet}$-module de Tor-dimension finie, modulo
les complexes essentiellement nuls il l'est et on peut identifier 
$\pi^* (\O_\la/\varpi_\la)$ au complexe 
$$
 [\underset{-1}{ \O_{\la\bullet}} \xrig{\;\times \varpi_\la\;} \underset{0}{
 \O_{\la\bullet}} ]
$$
puisque $\ker ( \O_{\la\bullet} \xrig{\;\times \varpi_\la\;}
 \O_{\la\bullet} )$ est essentiellement nul. 
Cela permet de définir pour $C\in \mathbb{D}^{e+}
(X^\N,{\O_{\la\bullet}})$ 
$$
\pi^* (\O_\la/\varpi_\la)
\overset{\mathbb{L}}{\otimes}_{\O_{\la\bullet}} C = \text{Tot}^\oplus
[ C\xrig{\; \varpi_\la\;} C]
$$
dans la catégorie $\mathbb{D}^{e+} (X^\N,{\O_{\la\bullet}})$ modulo
les complexes essentiellement nuls. On peut donc définir ce que cela
veut dire que ce complexe est essentiellement constant. On peut
également définir son image par $R\pi_*$. 

\subsection{$\O_\la$-complexes}

\begin{defi}[\cite{Eke} déf. 2.1]
Un $\O_{\la}$-complexe est  un $C\in \mathbb{D}^{e+} (
X^\N,\O_{\la\bullet})$ tel que $\pi^* (\O_\la/\varpi_\la)
\overset{\mathbb{L}}{\otimes}_{\O_{\la\bullet}} C$ soit
essentiellement constant 
\end{defi}

Les $\O_\la$-complexes sont les objets de base introduits par Ekedahl
afin de définir la catégorie dérivée $\la$-adique.

\subsection{La catégorie $\mathbb{D}^+ (X,\O_\la)$}

\begin{defi}[\cite{Eke} déf. 2.1] 
Un objet $C\in \mathbb{D}^{e+} (
X^\N,\O_{\la\bullet})$ est négligeable si  $\pi^* (\O_\la/\varpi_\la)
\overset{\mathbb{L}}{\otimes}_{\O_{\la\bullet}} C$ est essentiellement
nul
et un morphisme $f$ est essentiellement un isomorphisme si 
un cône de $f$ est négligeable.
\end{defi}

\begin{rema}
Les morphismes qui sont essentiellement des isomorphismes permettent
un calcul des fractions à gauche et à droite dans $ \mathbb{D}^{e+} (
X^\N,\O_{\la\bullet})$ et la catégorie localisée est triangulée. 
 En effet, le foncteur $C\mapsto \mathcal{H}^0
\left (  \pi^* (\O_\la/\varpi_\la)
\overset{\mathbb{L}}{\otimes}_{\O_{\la\bullet}} C \right )$ à valeurs
dans la catégorie des systèmes projectifs $\la$-adiques sur $X_\et$
modulo les complexes essentiellement nuls est un foncteur
cohomologique (cf. la remarque \ref{reth}).
\end{rema}

\begin{defi}[\cite{Eke} déf. 2.5]
La catégorie des $\O_\la$-complexes a pour objets les
$\O_\la$-complexes munis des morphismes 
$$
\Hom_{\O_\la-\text{complexes}} (A,B) = \Hom_{\text{pro}-\mathbb{D}^{e+}} (\tau (A),\tau(B))$$
Ekedahl pose alors $\mathbb{D}^+ (X,\O_\la) = $ la catégorie des
$\O_\la$-complexes localisée par rapport aux flèches qui sont 
essentiellement des isomorphismes. C'est une catégorie triangulée.
\end{defi} 

\subsection{Complexes normalisés}

Ekedahl introduit des représentants particuliers des $\O_\la$
complexes dans leur classes d'isomorphisme appelés complexes
normalisés. 

\begin{defi}[\cite{Eke} section 2] 
Pour $\F\in  \mathbb{D}^{e+} (
X^\N,\O_{\la\bullet})$ on note 
$$
\widehat{\F}= \mathbb{L}\pi^* R\pi_* \F\in  \mathbb{D}^{+} (
X^\N,\O_{\la\bullet})
$$
\end{defi}
Le foncteur $\F\mapsto \widehat{\F}$ se factorise par la catégorie localisée
par les morphismes qui sont essentiellement des isomorphismes. Cela
résulte des formules (\cite{Eke} lemme 1.5.ii)) 
$$
\forall n \;\; i_n^* \widehat{\F} = R\pi_* \left ( \pi^* (\O_\la/\varpi_\la^n)
\overset{\mathbb{L}}{\otimes}_{\O_{\la\bullet}} \F\right )
$$
où $i_n : \widetilde{X}_\et \ldrt \widetilde{X}_\et^\N$ est le
morphisme de topos ``étage $n$''. 
Ekedahl démontre alors qu'un complexe $\F$ est un $\O_\la$-complexe
ssi la flèche d'adjonction $\widehat{\F}\ldrt \F$ est essentiellement
un isomorphisme (\cite{Eke} proposition 2.2.i)).

\begin{defi}
Les complexes $\F\in \mathbb{D}^+(X^\N,\O_{\la\bullet})$ tels que
$\widehat{\F}\drt \F$ soit un isomorphisme dans $\mathbb{D}^+
(X^\N,\O_{\la\bullet})$ 
sont des $\O_\la$-complexes appelés complexes
normalisés. On note $\mathbb{D}^+(X^\N,\O_{\la\bullet})_{norm}$ la
sous-catégorie des complexes normalisés. 
\end{defi}

On renvoi à la proposition 2.2.ii) de \cite{Eke} pour une caractérisation
intrinsèque des complexes normalisés.

\begin{prop}[\cite{Eke}]
La catégorie $\mathbb{D}^+(X^\N,\O_{\la\bullet})_{norm}$ est l'image
essentielle du foncteur $\mathbb{L}\pi^* : \mathbb{D}^+
(X,\O_\la)_{litt}$ (la catégorie dérivée usuelle des $\O_\la$-modules
sur $X_\et$). De plus les foncteurs $\F\mapsto \widehat{\F}$ et
le foncteur canonique $\mathbb{D}^+(X^\N,\O_{\la\bullet})_{norm} \drt 
\mathbb{D}^+ (X,\O_\la)$ induisent une équivalence de catégories entre
$\mathbb{D}^+(X^\N,\O_{\la\bullet})_{norm}$ et  $\mathbb{D}^+ (X,\O_\la)$.
\end{prop}

\begin{rema}
Soit $\mathbb{D}^+ (X,\O_\la)_{litt}$ la catégorie dérivée ``usuelle''
des complexes de $\O_\la$-faisceaux sur $X_\et$. 
Le foncteur $\mathbb{L}\pi^*$ permet d'identifier $\mathbb{D}^+
(X,\O_\la)$ à la catégorie quotient de $\mathbb{D}^+
(X,\O_\la)_{litt}$ par les complexes $\F$ tels que
$\mathbb{L}\pi^*\F=0$.
\end{rema}

\subsection{La catégorie $\mathbb{D}^b (X,\O_\la)$}

\begin{defi}
On note $\mathbb{D}^b (X,\O_\la)$ la sous-catégorie de $\mathbb{D}^+
(X,\O_\la)$ formée des complexes dans $\mathbb{D}^{-}$. 
\end{defi}

L'un des points clef qui simplifie nettement la théorie par rapport au
cas de $\mathbb{D}^+$ est le (i) de la proposition 3.4 de \cite{Eke} : 
\begin{prop}[\cite{Eke}]
Un complexe $\F\in \mathbb{D}^{-} (X^\N,\O_{\la\bullet})$ est
négligeable ssi il est essentiellement nul. 
\end{prop}

\begin{coro}
La catégorie $\mathbb{D}^b (X,\O_\la)$ s'identifie au localisé de la
catégorie des $\O_\la$-complexes dans $\mathbb{D}^{-}$ relativement
aux flèches dont le cône est essentiellement nul. 
\end{coro}

\begin{lemm}\label{dokgon}
Soit $\F\in \mathbb{D}^b (X,\O_\la)$. Alors, $\widehat{\F}\in
\mathbb{D}^b (X^\N,\O_{\la\bullet})_{norm}$. Si de plus $\forall i
\notin [a \,b ]\;\;\mathcal{H}^i (\F)$ est essentiellement nul alors 
$\forall i\notin [a-1\; b]\; \mathcal{H}^i (\widehat{\F})$
est nul.
\end{lemm}
\dem
On utilise les formules 
$$
\forall n\;\; i_n^* \widehat{\F} = R\pi_* \left (  \pi^* (\O_\la/\varpi_\la^n)
\overset{\mathbb{L}}{\otimes}_{\O_{\la\bullet}} \F\right )
$$
$$
\pi^* (\O_\la/\varpi_\la^n)
\overset{\mathbb{L}}{\otimes}_{\O_{\la\bullet}} \F \simeq [\underset{-1}{\F}\xrig{\;\times\varpi_\la^n\;} \underset{0}{\F}]
$$
le dernier isomorphisme étant modulo les complexes essentiellement nuls.
Soit donc  $$f: \pi^{-1} C \drt [{\F}\xrig{\; \varpi_\la^n\;} {\F} ] $$
un morphisme de complexes de cône essentiellement nul. Alors, 
$$
i_n^*\widehat{\F} = C
$$
De plus, le cône de $f$ étant essentiellement nul on en déduit que
$\forall i\notin [a-1\; b] \; \mathcal{H}^i (\pi^{-1} C)$ est
essentiellement nul et donc 
$$
\forall i\notin [a-1\; b] \; \mathcal{H}^i (C)=0
$$
\qed

\begin{coro}
Le foncteur de normalisation induit une équivalence entre
$\mathbb{D}^b (X,\O_\la)$ et \\ $\mathbb{D}^b
(X^\N,\O_{\la\bullet})_{norm}$. 
\end{coro}

\subsection{La catégorie $\mathbb{D}^b_c (X,\O_\la)$ et sa t-structure tautologique}

Désormais on supposera toujours que $X$ est localement noethérien.

Rappelons qu'Ekedhal pose pour $\F\in \mathbb{D}^b (X,\O_\la)$
$$
\O_\la/\varpi_\la\overset{\mathbb{L}}{\otimes}_{\O_\la} \F := i_1^*\widehat{\F} = \O_\la/\varpi_\la \overset{\mathbb{L}}{\otimes}_{\O_\la} R\pi_* \F \in \mathbb{D}^b_c (X,\O_\la/\varpi_\la)
$$
\begin{defi}[\cite{Eke}] On pose
$$
\mathbb{D}^b_c (X,\O_\la) = \{\; \F\in \mathbb{D}^b (X,\O_\la) \;|\; \O_\la/\varpi_\la \overset{\mathbb{L}}{\otimes}_{\O_\la} \F\in \mathbb{D}^b_c (X,\O_\la/\varpi_\la) \; \}
$$
\end{defi}

\begin{defi}\label{kmpdgje}
\begin{itemize}
\item
Un faisceau $\la$-adique constructible sur $X_{\et}$ est un
$\O_{\la\bullet}$-module sur $X_{\et}$, $(\F_n)_n$ tel que $\forall
n\; \F_{n+1}\otimes \O_\la/\varpi_\la^n\iso \F_n$ 
et  $\F_1$ soit constructible. 
\item La catégorie des faisceau $\la$-adiques
constructibles est la sous-catégorie des $\O_{\la\bullet}$-modules
ayant pour objets les faisceaux $\la$-adique constructibles
\item Un faisceau essentiellement $\la$-adique constructible est un
  $\O_{\la\bullet}$-module isomorphe modulo les complexes
  essentiellement nuls à un faisceau   $\la$-adique constructible
\item La catégorie des faisceaux essentiellement $\la$-adiques
  constructibles est la sous-catégorie de la catégorie quotient des
  $\O_{\la\bullet}$-modules par les modules essentiellement nuls ayant
  pour objets les  faisceaux essentiellement $\la$-adique constructible
\end{itemize}
\end{defi}

Il est aisé de vérifier que le
 foncteur naturel de la catégorie des faisceaux $\la$-adiques
 constructibles vers la catégorie des faisceaux essentiellement
 $\la$-adiques constructibles induit une équivalence de catégories.

Rappelons le théorème principal de \cite{Jou1} :

\begin{theo}[\cite{Jou1}]
La catégorie des faisceaux (essentiellement) $\la$-adiques
constructibles est abélienne noethérienne et est équivalente à la
catégorie des faisceaux localement AR-$\la$-adiques constructibles.
\end{theo}

Ici la condition $AR$-nul (pour Artin-Rees) est plus forte que la
condition essentiellement nul. Si $X$ est quasicompact, si $f:\F \ldrt \mathcal{G}$ est un
morphisme de faisceau $\la$-adique constructibles et si $\mathcal{H}$
désigne le noyau de $\F$ dans la catégorie des
$\O_{\la\bullet}$-modules, celui-ci vérifie la condition de
Mittag-Lefler-Artin-Rees. Si $\mathcal{H}'$ désigne le système
projectif des images universelle de $\mathcal{H}$ alors 
 pour $r>>0$ le
$\O_{\la\bullet}$-module $\mathcal{H}[r]\otimes \O_{\la\bullet}$ (où l'on pose
$\mathcal{H}[r]_i=\mathcal{H}_{r+i}$) est
$\la$-adique constructible  égal au noyau de $f$ dans la catégorie
des faisceaux $\la$-adiques constructibles.

On renvoie plus généralement à \cite{Jou1} pour les propriétés de base des faisceaux $\la$-adiques constructibles.

\begin{prop}[\cite{Eke} section 3] \hspace{3mm}
 \begin{itemize}
\item Pour $(\F_n)_n$ un faisceau essentiellement $\la$-adique constructible, $(\F_n)_n$ placé en degré $0$ est un $\O_\la$-complexe.
\item Un $\F\in \mathbb{D}^b (X,\O_\la)$ appartient à $\mathbb{D}^b_c
  (X,\O_\la)$ ssi $\forall i \; \mathcal{H}^i (\F)$ est un faisceau  essentiellement $\la$-adique constructible.
\item Si $(\F_n)_n$ est un faisceau essentiellement $\la$-adique constructible et $\F$ le complexe concentré en degré $0$ associé
alors $\mathcal{H}^0 (\widehat{\F})$ est un faisceau $\la$-adique
constructible tel que les noyaux et conoyaux de l'application
$\mathcal{H}^0 (\widehat{\F}) \drt (\F_n)_n$ sont essentiellement
nuls. Le foncteur $\F\mapsto \mathcal{H}^0 (\widehat{\F})$ induit une
équivalence entre la catégorie des faisceaux essentiellement
$\la$-adiques constructibles  et la catégorie des faisceaux
$\la$-adiques constructibles, inverse de l'équivalence naturelle
décrite précédemment après la définition \ref{kmpdgje}.
\end{itemize}
\end{prop}

\begin{theo}
Soit $X$ un schéma localement noethérien  et $\mathcal{D}
=\mathbb{D}^b_c (X,\O_\la) $. Posons
$$
\mathcal{D}^{\leq 0} = \{ \;\F\in\mathcal{D}\;|\; \forall i>0\;
\mathcal{H}^i (\F)  \text{ est essentiellement nul }
 \}
$$
$$
\mathcal{D}^{\geq 0} = \{ \;\F\in\mathcal{D}\;|\;\forall i<0\;
\mathcal{H}^i (\F)  \text{ est essentiellement nul }
 \}
$$
Alors, $(\mathcal{D}^{\leq 0},\mathcal{D}^{\geq 0})$ est une
t-structure sur $\mathcal{D}$ de coeur la catégorie des faisceaux
essentiellement $\la$-adiques constructibles. Le foncteur
cohomologique associé est $\F\mapsto \mathcal{H}^0 (\F)$ et les
opérateurs de troncature $\tau_{\leq i}, \tau_{\geq i}$ sont ceux
induits par les opérateurs de troncatures usuels.

Si $\mathcal{D}_2 =\mathbb{D}_c^b( X,\O_{\la\bullet})_{norm}$, via
l'équivalence de catégories triangulées $\F\mapsto \widehat{\F}$ 
entre $\mathcal{D}$ et
$\mathcal{D}_2$, la t-structure induite sur $\mathcal{D}_2$ est
$$
\mathcal{D}_2^{\leq 0} = \{ \;\F\in\mathcal{D}_2\;|\ \forall i>0\;
\mathcal{H}^i (\F) =0\; \}
$$
$$
\mathcal{D}_2^{\geq 0} = \{ \;\F\in\mathcal{D}_2\;|\; \forall i<-1\;
\mathcal{H}^i (\F) =0 \text{ et } \mathcal{H}^{-1} (\F)
\text{ est ess. nul }\}
$$
De plus les opérateurs de troncature associés sont $\F\mapsto \tau_{\leq i} \F$ et $ \F\mapsto \widehat{\tau_{\geq i} \F}$. Le coeur de $\mathcal{D}_2$ s'identifie à la catégorie des faisceaux $\la$-adiques constructibles et le foncteur cohomologique associé est $\F\mapsto \mathcal{H}^0 (\widehat{\mathcal{H}^0 (\F)})$. 
\end{theo}
\dem
C'est une conséquence du théorème précédent couplé au lemme \ref{dokgon}.
\qed

\subsection{La t-structure perverse sur $\mathbb{D}^b_c (X,\O_\la)$}

\begin{defi}
Soit $\mathcal{D}=\mathbb{D}^b_c (X,\O_\la)$. Posons 
$$
\,^p\mathcal{D}^{\leq 0} = \{ \;\F\in\mathcal{D}\;|\; \forall x\in X\; \forall i>-\dim \overline{\{x \}}\; \mathcal{H}^i (i_x^* \F) \text{ est ess. nul }\}
$$
$$
\,^p\mathcal{D}^{\geq 0} = \{ \;\F\in\mathcal{D}\;|\; \forall x\in X\; \forall i<-\dim \overline{\{x \}}\; \mathcal{H}^i (R i_x^! \F) \text{ est ess. nul }\}
$$
\end{defi}

\begin{theo}\label{experv}
Supposons $X$ noethérien et que $\mathbb{D}^b_c (X,\O_\la/\varpi_\la)$ possède un module dualisant. Alors $(\,^p\mathcal{D}^{\leq 0},\,^p\mathcal{D}^{\geq 0}) $ est une t-structure sur $\mathcal{D}$. Son coeur noté $\text{Perv} (X,\O_\la)$ est une catégorie abélienne $\O_\la$-linéaire noethérienne. 

De plus, si $\mathcal{D}_2 = \mathbb{D}^b_c (X^\N,\O_{\la\bullet})_{norm}$, dans l'équivalence $\F\mapsto \widehat{\F}$ entre $\mathcal{D}$ et $\mathcal{D}_2$ la t-structure précédente correspond à 
$$
\,^p\mathcal{D}_2^{\leq 0} =  \{ \;\F\in\mathcal{D}_2\;|\; \forall x\in X\; \forall i>-\dim \overline{\{x \}}\; \mathcal{H}^i (i_x^* \F) =0 \;\}
$$
$$
\,^p\mathcal{D}_2^{\geq 0} = \{ \;\F\in\mathcal{D}_2\;|\; \forall x\in X\; \forall i<-\dim \overline{\{x \}}-1 \; \mathcal{H}^i (R i_x^! \F) = 0 \text{ et } \mathcal{H}^{-\dim \overline{\{x \}}-1} (R i_x^! \F ) \text{ est ess. nul } \}
$$
\end{theo}
\dem
Montrons que $(\,^p\mathcal{D}^{\leq 0},\,^p\mathcal{D}^{\geq 0})$ est une t-structure par récurrence noethérienne. Supposons donc que pour tout sous-schéma fermé $F\subsetneq X$ la définition donnée définisse une t-structure sur $\mathbb{D}^b_c(F,\O_\la)$ (noter qu'un tel schéma $F$ vérifie bien les hypothèses de l'énoncé puisque si $i:F\hookrightarrow X$ et $K_X$ est dualisant sur $X$ alors $Ri^! K_X$ est dualisant sur $F$). 
 Tout $\F\in\mathcal{D}$ vérifie qu'il existe un ouvert dense $U$ dans $X$ tel que $\forall i \; \mathcal{H}^i (\F)_{|U}$ soit un système local $\la$-adique modulo les systèmes essentiellement nuls (\cite{Jou1} proposition 1.2.6). Soit $K_X\in \mathbb{D}^b_c (X,\O_\la/\varpi_\la)$ un module dualisant. Il existe un ouvert dense $V$ dans $X$ et un entier $\delta\in \Z$ tels que $K_{X|V}\simeq \L [\delta]$ où $\L$ est un $\O_\la/\varpi_\la$-module localement libre de rang $1$.
\\
 Pour un ouvert $U$ irréductible dans $X$ tel que le module dualisant $K_{X|U}$ soit de la forme précédente on vérifie facilement que si $\mathcal{D'}_U$ est la sous-catégorie de $\mathbb{D}^b_c (U,\O_\la)$ définie par $\F\in \mathcal{D}'_U$ ssi $\forall i\;\mathcal{H}^i (\F)$ est essentiellement un système local $\la$-adique alors la t-structure définie dans l'énoncé définit une t-structure sur $\mathcal{D}'_U$. On utilise pour cela le fait que sur $U$ le module $\O_\la/\varpi_\la\in\mathbb{D}^b_c (U,\O_\la/\varpi_\la)$ est dualisant et que donc grâce à la formule (cf. le chapitre 4 de \cite{Eke}) 
 $$\O_\la/\varpi_\la\overset{\mathbb{L}}{\otimes}_{\O_\la} R\underline{\Hom}_{\O_\la} (\F,\G)
\simeq  R\underline{\Hom}_{\O_\la/\varpi_\la}(\O_\la/\varpi_\la\overset{\mathbb{L}}{\otimes}_{\Zl} \F,\O_\la/\varpi_\la\overset{\mathbb{L}}{\otimes}_{\Zl}\G) 
$$
$\O_{\la\bullet} \in \mathbb{D}^b_c (X,\O_\la)$ est dualisant en un sens évident et échange les foncteurs $i_x^*$ et $Ri_x^!$. 
Sur $\mathcal{D}'_U$ les opérateurs de troncature pervers sont alors les opérateurs de troncature usuels décalés par la dimension de $U$).
\\
Pour un ouvert $U$ comme précédemment notons $\widetilde{\mathcal{D}}
(U)$ la sous-catégorie triangulée de $\mathcal{D}$ formée des
$\F\in\mathcal{D}$ tels que $\forall i\; \mathcal{H}^i(\F)_{| U}$ soit essentiellement un système local $\la$-adique. 
Si $j:U\hookrightarrow X$ et $i:X\setminus U \hookrightarrow X$ alors, l'existence du module dualisant implique que les foncteurs $Ri^!$ et $Rj_*$ respectent les catégories $\mathbb{D}^b_c$ pour les coefficients de torsion et donc pour les coefficients $\la$-adiques. 
Alors, la définition donnée dans l'énoncé induit une t-structure sur $\widetilde{\mathcal{D}} (U)$ par recollement (cf. \cite{BBD}) de la t-structure sur $\mathcal{D}'_U$ et celle sur $\mathbb{D}^b_c (X\setminus U,\O_\la)$ obtenue  par hypothèse de récurrence. Étant donné que ``$\mathcal{D}=\bigcup_U \widetilde{\mathcal{D}} (U)$'' on en déduit facilement que l'on a bien une t-structure sur $\mathcal{D}$.
\\

La description de la t-structure sur $\mathbb{D}^b_c
(X^\N,\O_{\la\bullet})_{norm}$ résulte de ce que $\forall \F \in
\mathbb{D}^b(X,\O_\la)\;$ $ \forall x\in X\; i_x^* (\widehat{\F}) = \widehat{i_x^*\F}$, $
  Ri_x^! (\widehat{\F}) = \widehat{Ri_x^!\F}$ et du lemme \ref{dokgon}.

La noethérianité de la catégorie $\text{Perv} (X,\O_\la)$ se démontre
comme pour le théorème 8.3 de \cite{Gabber1} en utilisant la
noethérianité de la catégorie des faisceaux $\la$-adiques constructibles.
\qed

\subsection{La t-structure perverse sur $\mathbb{D}^b(X^\N,\O_{\la\bullet})$}

Définissons maintenant une t-structure ``naïve'' sur
$\mathbb{D}^b(X^\N,\O_{\la\bullet})$. 

\begin{theo}\label{roghonut}
Soit $X$ un schéma noethérien de dimension finie. Soit
$\mathcal{D}=\mathbb{D}^b (X^\N,\O_{\la\bullet})$. Posons 
$$
\,^p\mathbb{D}^{\leq 0} = \{\;\F\in \mathbb{D}^b(X^\N,\O_{\la\bullet})
\;|\; \forall x\in X\;\forall i>-\dim \overline{\{ x\} }\;
\mathcal{H}^i (i_x^* \F) = 0 \;\}
$$
$$
\,^p\mathbb{D}^{\geq 0} = \{\;\F\in \mathbb{D}^b(X^\N,\O_{\la\bullet})
\;|\; \forall x\in X\;\forall i<-\dim \overline{\{ x\} }\;
\mathcal{H}^i (Ri_x^! \F) = 0\;\} 
$$
Alors, $(\,^p\mathbb{D}^{\leq 0}, \,^p\mathbb{D}^{\geq 0} )$ forme une
t-structure sur $\mathcal{D}$.
\end{theo}
\dem
Les méthodes de l'article \cite{Gabber1} s'adaptent aussitôt.
\qed

\begin{defi}
On note $\text{Perv} ( X^\N,\O_{\la\bullet})$ le coeur de la
t-structure précédente. C'est une catégorie abélienne $\O_\la$-linéaire.
\end{defi}

\subsection{Lien entre les deux t-structures}

Dans cette section le schéma $X$ satisfait aux hypothèses du théorème
\ref{experv}.

\begin{prop}
Soit $\F\in \text{Perv} (X,\O_\la)$ sans $\varpi_\la$-torsion. Alors 
$\widehat{\F}\in \text{Perv} (X^\N,\O_{\la\bullet})$.
\end{prop}
\dem
D'après la description de la t-structure sur les faisceaux normalisés 
donnée dans le théorème 
 \ref{experv} $\;\widehat{\F}\in \,^p\left
  (\mathbb{D}^b (X^\N,\O_{\la\bullet}) \right )^{\leq 0}$.
Soit $\mathcal{G}$ un cône de l'application $\widehat{\F}\xrig{\times
  \varpi_\la} \widehat{\F}$. Par hypothèse 
$$\mathcal{G}\in \text{Perv}
(X,\O_\la) \cap \mathbb{D}^b (X^\N,\O_\la)_{norm}$$
Soit  $x\in X$ et $i<-\dim \overline{\{ x\} }$. Si $i<-\dim
\overline{\{ x\} } -1$ d'après le théorème \ref{experv} $\mathcal{H}^i
( Ri_x^! \widehat{\F})=0$. Il y a de plus une suite exacte 
$$\underbrace{
\mathcal{H}^{-\dim \overline{\{ x\} }-2} ( Ri_x^!\mathcal{G})}_0 \ldrt 
\mathcal{H}^{-\dim \overline{\{ x\} }-1} ( Ri_x^! \F) \xrig{\;\times
  \varpi_\la\;} \mathcal{H}^{-\dim \overline{\{ x\} }-1} ( Ri_x^! \F)
$$
et donc la multiplication par $\varpi_\la$  sur le $\O_{\la\bullet}$-module 
$\mathcal{H}^{-\dim \overline{\{ x\} }-1} ( Ri_x^! \F)$ est injective,
module 
qui est donc
nul.
\qed

\begin{coro}\label{cortyu}
Le foncteur $\F\mapsto \widehat{\F}$ de $\mathbb{D}^b_c
(X,\O_\la)\otimes L_\la$
dans $\mathbb{D}^b (X^\N,\O_{\la\bullet})\otimes L_\la$ est t-exact et
induit donc un foncteur exact 
$\text{Perv} (X,\O_\la) \otimes L_\la \ldrt \text{Perv}
(X^\N,\O_{\la\bullet})\otimes L_\la$. 
\end{coro}
\dem
Si $\F\in \text{Perv} (X,\O_\la)$, $\F$ étant un objet noethérien il 
existe un entier $N\in \N$ tel que $\F/\F[\varpi_\la^N]$ soit sans
$\varpi_\la$-torsion. L'existence du foncteur exacte entre 
$\text{Perv} (X,\O_\la) \otimes L_\la$ et  $\text{Perv}
(X^\N,\O_{\la\bullet})\otimes L_\la$ en résulte. La première assertion
sur la t-exactitude s'en déduit par dévissage en ``découpant un objet en
faisceaux pervers décalés par troncatures succesives''. 
\qed

\subsection{Catégories dérivées filtrées $\la$-adiques}

\begin{defi}
On note $\mathbb{D} F (X^\N,\O_{\la\bullet})$ la catégorie dérivée
filtrée des complexes de $\O_{\la\bullet}$-modules sur $X_\et^\N
= ( \dots \drt X_\et \drt \dots \drt X_\et )$
où les filtrations sont prises décroissantes finies (et non pas finies degré par degré). On
note
$\mathbb{D}^{e+} F (X^\N,\O_{\la\bullet})$ la sous-catégorie formée 
des complexes filtrés $C$ tels que pour $i<<0$  $\;\forall j\;$
$\;\mathcal{H}^i ( Gr^j C)$ soit essentiellement nul. 
\end{defi}

Il y a un couple de foncteurs 
$$
\xymatrix{
\mathbb{D}^{e+} F ( X^\N,\O_{\la\bullet}) \ar@<1ex>[r]^{R\pi_*} 
&  \ar@<1ex>[l]^{\mathbb{L}\pi^*}  \mathbb{D}^+ F( X,\O_{\la})_{litt}
}
$$
\\
Tout objet de $\mathbb{D^+}F (X^\N,\O_{\la\bullet})$ est isomorphe à
un complexe $\Fil^\bullet (I^\bullet_n)_n$ où $\forall i\forall j\;
Gr^i (I^{j}_n)$ est un faisceau injectif et $\forall n\;\;
Gr^i (I^j_{n+1})\drt Gr^i (I^j_n )$ est un épimorphisme qui possède 
une section. Pour un tel complexe $R\pi_* ( \Fil^\bullet
I^\bullet_\bullet)
= \underset{n}{\limp} \Fil^\bullet I^\bullet_n$. On définit alors 
$R\pi_*$ sur $\mathbb{D}^{e+} (X^\N,\O_{\la\bullet})$ par une
procédure analogue à celle de \cite{Eke} en utilisant le pro-système 
$\tau (-)$. 
\\
Quant au foncteur $\mathbb{L}\pi^*$ il se calcule naturellement sur
les complexes filtrés $\Fil^\bullet C^\bullet$ tels que 
$\forall i,j\; Gr^i C^j$ soit $\O_{\la}$-plat ( ce qui est possible
puisque $\O_{\la\bullet}$ est un $\pi^{-1}(\O_\la)$-module de
Tor-dimension finie).

\begin{defi}
Un complexe filtré $C\in \mathbb{D}^{e+}F ( X^\N,Ab)$ est
essentiellement nul si $\forall j\; Gr^j C$ l'est.
 Il est
essentiellement constant s'il est isomorphe dans la catégorie 
des complexes filtrés sur $X_\et^\N$ localisée par les flèches de cône
essentiellement nul à un objet de la forme $\pi^{-1} D$ où $D$ est un
complexe filtré de groupes abéliens sur $X_\et$. 
\end{defi}

\begin{defi}
\begin{itemize}
\item
Un $\O_\la$-complexe filtré est un $C\in\mathbb{D}^{e+} F
(X^\N,\O_{\la\bullet})$ tel que  le complexe filtré $\pi^* (\O_{\la}/\varpi_\la)
 \overset{\mathbb{L}}{\otimes}_{\O_\la\bullet} C$ soit essentiellement
 constant.
\item Un complexe filtré normalisé est un $C$ tel que  
le morphisme 
$\tau ( \mathbb{L}\pi^* R\pi_* C) \drt \tau (C)$ soit un
isomorphisme de pro-objets dans $\mathbb{D}^+ F
(X^\N,\O_{\la\bullet})$. 
\item Le complexe filtré $C$ est négligeable si $\pi^*
  (\O_\la/\varpi_\la)
 \overset{\mathbb{L}}{\otimes}_{\O_\la\bullet } C$ est essentiellement
  nul.
\item Un morphisme est essentiellement  un isomorphisme  si son cône
  est  négligeable. 
\item Un morphisme de $\O_\la$-complexe est un morphisme entre
  pro-systèmes $\tau(-)$ associés dans $\mathbb{D}^{e+}F(X^\N,\O_{\la\bullet})$.
\end{itemize}
\end{defi}

\begin{defi}
On note $\mathbb{D}^+ F (X,\O_\la)$ le localisé de la catégorie 
des $\O_\la$-complexes vis à vis des morphismes qui sont
essentiellement
des isomorphismes.
\end{defi}

Comme dans \cite{Eke}, la sous-catégorie des complexes normalisés dans 
$\mathbb{D}^+F (X^\N,\O_\la)$ s'identifie à $\mathbb{D}^+ F (X,\O_\la)$.
Cette catégorie est munie de foncteurs
$$
\forall j \; Gr^j : \mathbb{D}^+F (X,\O_\la) \ldrt  \mathbb{D}^+ (X,\O_\la)
$$

\begin{prop} \label{olpic}
Reprenons les hypothèses du théorème \ref{experv}. 
Il y a une équivalence de catégories entre la catégorie des objets de
$\text{Perv} (X,\O_\la)$ munis d'une filtration finie et celle des 
objets $C$ de $\mathbb{D}^+ F (X,\O_\la)$ tels que $\forall i\; Gr^iC\in
\text{Perv} (X,\O_\la)$. 
\end{prop}
\dem
La démonstration de la proposition \ref{PervFilt} s'adapte. Pour la
pleine fidélité il suffit essentiellement de connaître l'existence d'une suite
spectrale 
$$
\forall A,B\in \mathbb{D}^+ F (X,\O_\la)\;\; E^{pq}_1 \limpl
\Hom_{\mathbb{D}^+F (X,\O_\la)} (A, B[p+q])
$$
où
$$
E^{pq}_1 = \left \{ \dpt{\prod_{j-i=p} \Hom_{\mathbb{D}^+ (X,\O_\la)} ( Gr^i A, Gr^j B[p+q])
\;\text{ si } p \geq 0} \atop
 0\; \text{ si } p<0 \right.
$$
Mais si l'on note $\forall C\in \mathbb{D}^+ F (X,\O_\la)\;\;
\widehat{C} = \mathbb{L}\pi^* R\pi_* C\in \mathbb{D}^+
(X^\N,\O_{\la\bullet})$ alors 
$$
\Hom_{\mathbb{D}^+F(X,\O_\la)} (A,B) = \Hom_{\mathbb{D}^+F
  (X^\N,\O_{\la\bullet})} (\widehat{A},\widehat{B})
$$
L'existence de la suite spectrale se déduit alors de celle d'Illusie
pour la catégorie dérivée filtrée $\mathbb{D}^+ F
(X^\N,\O_{\la\bullet})$ une fois vérifié que 
$$
\forall C\in \mathbb{D}^+F (X,\O_\la)\;\; Gr^i \widehat{C}=
\widehat{Gr^i C}
$$
La surjectivité essentielle se démontre exactement de la même façon
que dans la proposition \ref{PervFilt}.
\qed

On vérifie maintenant aussitôt le lemme suivant :

\begin{lemm}
Soit $\mathcal{A}$ une catégorie abélienne $\O_\la$-linéaire et
$F\mathcal{A}$ la catégorie formée des objets de $\mathcal{A}$ munis
d'une filtration de longueur finie. Le foncteur
$$
\left ( F \mathcal{A}\right ) \otimes_{\O_\la} L_\la \ldrt F \left (
  \mathcal{A}\otimes_{\O_\la} L_\la \right )
$$
induit une équivalence de catégories.
\end{lemm}

\begin{defi}
Pour $\mathcal{A}$ une catégorie abélienne $\O_\la$-linéaire on note
indifféremment 
$F\mathcal{A} \otimes L_\la$ l'une des deux catégories définies dans
le lemme précédent et identifiées grâce à celui-ci. 
\end{defi}

La vérification du lemme qui suit est également immédiate.

\begin{lemm}\label{ducasse}
Soit $X$ un schéma noethérien de dimension finie tel que
$\mathbb{D}^b_c (X,\O_\la/\varpi_\la)$ possède un module dualisant. Le
diagramme suivant est commutatif
$$
\xymatrix{
F\text{Perv} (X,\O_\la) \otimes L_\la \ar[r]^\a \ar[d]^\gamma & F\text{Perv}
(X^\N,\O_{\la\bullet}) \otimes L_\la \ar[d]^\delta \\
\mathbb{D}^b_c F (X,\O_\la)\otimes L_\la \ar[r]^\beta & \mathbb{D}^b F
(X^\N,\O_{\la\bullet})\otimes L_\la
}
$$ 
où $\a$ est le foncteur définit grâce au corollaire \ref{cortyu},
$\beta$ est le foncteur de normalisation $\F\mapsto \widehat{\F}$,
$\gamma$ est définit via la proposition \ref{olpic} et $\delta$ grâce
à la proposition \ref{PervFilt}. 
\end{lemm}

\section{Le théorème clef}

\begin{theo}\label{cleft}
Soit $X_0$ un $k$-schéma de type fini et $x\in X_0$ (un point
non nécessairement fermé). Soit $$
f : \spec (\widehat{\O}_{X_0,x}) \ldrt \spec ( \O_{X_0,x})
$$
Alors, si $\La$ est de torsion première à $p$ 
$$
f^* : \mathbb{D}^b_c ( \spec ( \O_{X_0,x})_{\text{ét}}, \La)
\ldrt \mathbb{D}^+ ( \spec (\widehat{\O}_{X_0,x})_{\text{ét}},\La)
$$
et
$$ 
f^* : \mathbb{D}^b ( \spec (
\O_{X_0,x})_{\text{ét}}^\N,\O_{\la\bullet}) \ldrt 
\mathbb{D}^+ ( \spec (\widehat{\O}_{X_0,x})_{\text{ét}}^\N,\O_{\la\bullet})
$$
sont t-exacts, où dans le second cas la t-structure est celle définie
dans le théorème \ref{roghonut}.
\end{theo}
\dem
Il est clair que le cas de torsion implique le second cas
puisqu'il suffit de le vérifier étage par étage. 
 Nous nous y
restreignons donc.
Notons 
$$
A= \O_{X_0,x}, \;\; Y=\spec (\widehat{A}), \;\; Z = \spec (A)
$$
 On vérifie facilement que 
$$
\forall y\in Y\;\; \dim \overline{\{f(y)\}} \geq \dim \overline{\{y\}}
$$
et que donc $f^*$ est t-exact à droite :
$$
f^* \left (\,^{p}\mathbb{D}^{\leq 0} \right ) \subset
\,^{p}\mathbb{D}^{\leq 0}
$$
La t-exactitude à gauche de $f^*$ repose sur le théorème de 
désingularisation de Popescu (\cite{Pop}). D'après celui-ci, 
l'anneau $A$ étant excellent, il existe un système inductif 
filtrant $(B_i)_{i\in I}$ de $A$-algèbres lisses et un isomorphisme 
de $A$-algèbres
$$
\widehat{A}\simeq \underset{i\in I}{\limi} B_i
$$
Quitte à localiser  les  $B_i$ on peut supposer que $\forall i\; B_i$
est un anneau local, les morphismes 
$$
A\ldrt B_i \ldrt \widehat{A}
$$
sont locaux et $B_i$ est une $A$-algèbre essentiellement lisse. Notons
$W_i = \spec (B_i)$ et 
$$
\xymatrix{
Y \ar[rr]^f \ar[rd]^{h_i} && Z  \\
 & W_i \ar[ru]^{g_i}
}
$$
Notons $d_i$ la dimension relative de $g_i $. Soit donc 
$$
\F\in \mathbb{D}^b_c (Z,\La) \cap \,^{p}\mathbb{D}^{\geq 0}
$$

Il existe un morphisme lisse de $k$-schémas de type fini de dimension
relative $d_i$ 
$$
\widetilde{W}_i \xrig{\; \widetilde{g}_i\;} \widetilde{Z}
$$
des points $a\in \widetilde{W}_i$ et $b\in \widetilde{Z}$
tels que $\widetilde{g}_i (a ) =b$, un 
complexe de faisceaux $\widetilde{\F}\in \mathbb{D}^b_c
(\widetilde{Z},\La) \cap \,^p \mathbb{D}^{\geq -\dim \overline{\{x\}}}$ et
un diagramme commutatif 
$$
\xymatrix{
W_i \ar[r]^{g_i}\ar[d]^\a_\simeq & Z \ar[d]^\beta_{\simeq}  \\
\spec ( \O_{\widetilde{W}_i,a}) \ar[r] \ar[d]^\gamma 
 & \spec (\O_{\widetilde{Z},b}) \ar[d]^\delta
\\
\widetilde{W}_i \ar[r]^{\widetilde{g}_i} & \widetilde{Z} 
}
$$ 
tel que 
$$
(\delta\beta )^* \widetilde{\F} = \F
$$
Le décalage $\dim \overline{\{x\}}= \dim \overline{\{b\}}$ provient de ce que
pour $U$ un schéma de type fini sur $k$ et $u\in U$ si 
$\ph : \spec (\O_{U,u}) \drt U$ alors $\ph^* [\dim \overline{ \{u \}}]$
est t-exact.
\\
D'après  \cite{BBD} pages 108-109, $\widetilde{g}_i$ étant lisse $\widetilde{g}_i^* [d_i]$ est t-exact et donc 
$$
\widetilde{g}_i^* \widetilde{\F} [d_i] \in \,^p\mathbb{D}^{\geq-\dim
  \overline{\{x \}}}
$$
 Étant donné que l'extension de corps résiduel 
induite par le morphisme d'anneaux locaux $A \drt B_i$ est triviale,
$\dim\overline{\{x \}} = \dim \overline{\{ a \}}$. On en déduit que
$$ \forall i \;\;
g_i^*\F [d_i] = \a^* \gamma^* \widetilde{g}_i^* \widetilde{F}[d_i]
\in \,^p\mathbb{D}^{\geq 0}
$$

Soit maintenant $\mathfrak{q}\in Y$. Il existe un $i_0\in I$ tel que
$$
\forall i\geq i_0 \text{ si } \mathfrak{p}_i = h_i (\mathfrak{q})
\text{  alors } \mathfrak{p}_i. \widehat{A}= \mathfrak{q}
$$
On vérifie alors que si
\begin{eqnarray*}
i_{\mathfrak{q}} &:& \spec ( k (\mathfrak{q})) \ldrt  Y \\
\mu_i &:& \spec (k (\mathfrak{q})) \ldrt \spec ( k (\mathfrak{p}_i))
\\
i_{\mathfrak{p}_i} &:& \spec ( k (\mathfrak{p}_i)) \ldrt W_i
\end{eqnarray*}
alors  
$$\boxed{
Ri_{\mathfrak{q}}^! \left ( f^*\F \right )
= \underset{i\geq i_0}{\limi} \mu_i^* R i_{\mathfrak{p}_i}^! \left (
g_i^*\F \right )}
$$
Étant donné que $g_i^*\F\in \,^p\mathbb{D}^{\geq - d_i}$ on en déduit 
que 
$$
\mathcal{H}^k \left ( R i_{\mathfrak{q}}^! (f^*\F) \right ) = 0
\;\text{  si  }  \forall i\geq i_0
\;\; k<-\dim \overline{\{\mathfrak{p}_i\}} + d_i
$$

Le théorème résulte donc de l'inégalité suivante 
$$ \forall i\geq i_0\;\;
\dim \overline{\{ \mathfrak{p}_i\}} \leq \dim
\overline{\{\mathfrak{q}\}} + d_i
$$

Celle-ci se démontre de la façon suivante : 
étant donné que l'extension de corps résiduels induite par le
morphisme local essentiellement lisse $A\drt B_i$
 est triviale, il y a un
diagramme commutatif de $A$-algèbres
$$
\xymatrix{
B_i \ar[r]\ar@{^(->}[d] & \widehat{A} \\
\widehat{B}_i \ar@{-}[r]^(0.3){\simeq} & \widehat{A} [[T_1,\dots, T_{d_i} ]]
\ar[u]^\psi
}
$$
où $\psi$ est continu, $\psi (T_i)$ appartenant à l'idéal
maximal de $\widehat{A}$. De plus c'est un morphisme de
$\widehat{A}$-algèbres.
Dès lors, si $I = \mathfrak{p}_i \widehat{B}_i \subset \widehat{A}
[[T_1,\dots, T_{d_i} ]]$, si $\forall k\; \psi (T_k)= t_k$ alors
$$
\ker \psi = (T_1-t_1,\dots, T_{d_i} - t_{d_i} )
$$
ce qui implique que
$$
\widehat{A}/\mathfrak{q} = \widehat{A}/\psi (I) \widehat{A}
\simeq \widehat{A} [[T_1,\dots, T_{d_i} ]] / \left ( I + 
(T_1-t_1,\dots,T_{d_i}-t_{d_i} )\right )
$$
et donc
$$
\dim \widehat{A}/\mathfrak{q} \geq \dim \underbrace{\widehat{A}
  [[T_1,\dots,T_{d_i} ]] /I}_{\simeq \widehat{B}_i /I \simeq
  \widehat{(B_i /\mathfrak{p}_i)}} -d_i = \dim B_i/\mathfrak{p}_i - d_i
$$
La première inégalité résultant de ce que 
$$
\text{ht}_{\widehat{A} [[T_1,\dots, T_{d_i} ]]/I} \left ( 
I+ (T_1-t_1,\dots,T_{d_i}-t_{d_i} )/I\right ) \leq d_i
$$
et qu'étant donné que $B_i$ est excellent $\widehat{B}_i/I= \widehat{B_i/\mathfrak{p}_i}$ est
équidimensionnel (EGA IV Scholie 7.8.3. (x)).
\qed

\section{Démonstration du théorème principal}

\subsection{Hypothèses et notations}

Soit $A$ une $\O$-algèbre plate locale et complète. On suppose que le
corps résiduel de $\O$ est séparablement clos i.e. $s=\overline{s}$ (on vérifie
aisément que l'on peut se ramener à ce cas là pour nos besoins). Soit 
$$
\F\in \mathbb{D}^b (\spec
(A[\frac{1}{\pi}])_{\et}^\N,\O_{\la\bullet})\otimes L_\la
$$
On suppose qu'il existe un $\O$-schéma de type fini $X$, un point
fermé $x\in X_s$ tel que $A\simeq \widehat{\O}_{X,x}$ et un faisceau
pervers 
$\G\in \text{Perv} (X_\eta,\O_\la)\otimes L_\la$ tel que via le
foncteur $\G\mapsto \widehat{\G} \in \Perv
(X_\eta^\N,\O_{\la\bullet})\otimes L_\la$ du corollaire \ref{cortyu}
et un isomorphisme $A\simeq
\widehat{\O}_{X,x}$ on ait
$$
\widehat{\G}_{|\spec (A[\frac{1}{\pi}])} \simeq \F
$$
On note encore $\F$ pour l'extension de $\F$ à $\spec (A)_{\bar{\eta}}$. 
Considérons le diagramme 
$$
\xymatrix{
\spec (A)_{\bar{\eta}} \ar@{^(->}[r]^j & \spec (A\otimes_\O
\overline{\O}) & \spec (A)_{s} \ar@{_(->}[l]_i
}
$$
où $\spec (A)_{s}= \spec (A/\pi A)$.

\subsection{Perversité des cycles évanescents formels et quasi-unipotence
  de la monodromie} 

\begin{theo}\label{exfup}
Soit $$\Rpsi (\F) = i^* R j_* \F\in  \mathbb{D}^+ ( \spec
(A)_{s}^\N,\O_{\la\bullet})\otimes L_\la$$
 Alors, 
$$
\Rpsi (\F)\in \text{Perv} (\spec
(A)_{s}^\N,\O_{\la\bullet})\otimes L_\la
$$
De plus, pour tout  $T\in I_{\bar{\eta}/\eta}$  l'action de $T$ sur
$\Rpsi (\F)$ est quasi-unipotente : pour $N$ suffisamment divisible et
$k>>0$ 
$$
(T^N - \text{Id})^k = 0 \in \End \left (\Rpsi (\F)\right )
$$
\end{theo}
\dem
Fixons un quadruplet $(X,x,\G,f)$ où $(X,x,\G)$ est comme
dans la section d'auparavant et
$$
f: \spec (A) \ldrt X
$$
se factorise par $\spec (\O_{X,x})\ldrt X$ et induit un isomorphisme 
$$
\spec (A) \iso \spec (\O_{X,x})
$$
D'après le théorème de changement de base régulier de Fujiwara
(corollaire 7.1.6 de \cite{Fuj2} où il est énoncé dans le cas de
torsion mais le cas utilisé ici s'en déduit aussitôt puisqu'il suffit
de le vérifier étage par étage dans le topos $X^{\N}_{\et}$) il y a un isomorphisme 
$$
\bar{f}^* \Rpsi ( \widehat{\G}) \simeq \Rpsi (\F)
$$
où $\bar{f}: \spec (A\otimes \bar{\O}) \ldrt X\otimes \overline{\O}$.
Le résultat se déduit alors du cas algébrique connu pour $X$
(i.e. $\Rpsi (\F)$ est de longueur finie, $\End (\Rpsi (\F))$ est un
$L_\la$-e.v. de dimension finie et l'action de l'inertie
$I_{\bar{\eta}|\eta}$ dessus est continue) puisque 
$$
\Rpsi (\widehat{\G}) = \widehat{\Rpsi (\G )}
$$
et 
du théorème \ref{cleft}. 
\qed

\subsection{\'Enoncé et démonstration du théorème principal}

Rappelons maintenant le lemme suivant :

\begin{lemm}
Soit $\mathcal{A}$ une catégorie abélienne $L_\la$-linéaire et $A\in
\mathcal{A}$. Soit $N\in \End (\mathcal{A})$ un endomorphisme
nilpotent. Il existe alors une unique filtration croissante finie
$\Fil_\bullet A$ sur $A$ telle que $\forall i\; N \Fil_i A \subset
\Fil_{i-2} A$ et $N$ induise un isomorphisme $N^i : \text{Gr}_i A \iso
\text{Gr}_{-i} A$. De plus, si $F: \mathcal{A} \ldrt \mathcal{B}$ est
un foncteur exacte entre catégories abéliennes alors l'objet filtré
associé à $(F(A),F(N))$ est $F(\Fil_\bullet A)$.
\end{lemm}

\begin{defi}
Soit $T\in I_{\bar{\eta}|\eta}$ s'envoyant sur un générateur de la
composante $\Z_\ell (1)$ dans l'inertie modérée. Soit $N=\log T\in
\End (\Rpsi (\F))$ . Soit $\Fil_\bullet \Rpsi (\F)$ l'objet filtré
associé par le lemme précédent. Nous noterons $K(\F)\in \mathbb{D}^b F
(\spec (A)_s^\N,\O_{\la\bullet})\otimes L_\la$ l'objet associé par la
proposition \ref{PervFilt} à $(\Fil_{-i} \Rpsi (\F))_{i\in \Z}$.
\end{defi}

\begin{prop}\label{kitulo}
Soit $(X,x,\G,f)$ comme dans la démonstration du théorème
\ref{exfup}. Considérons le complexe filtré  $K(\G)\in \mathbb{D}^b_c F (X,\O_\la)\otimes
L_\la$  associé à $\Rpsi (\G)$ muni de sa
filtration de monodromie par le lemme précédent et la proposition
\ref{olpic}.
Il y a un isomorphisme canonique 
$$
f^*\widehat{ K (\G)} \simeq K(\F)
$$
\end{prop}
\dem
La démonstration résulte du théorème \ref{cleft} couplé à
l'isomorphisme de Fujiwara comme dans la démonstration du théorème
\ref{exfup} ainsi que du lemme \ref{ducasse}. Il faut par exemple
vérifier que le théorème  \ref{cleft} implique la commutativité du
diagramme suivant 
$$
\xymatrix{
F\text{Perv} (X^\N,\O_{\la\bullet}) \otimes L_\la \ar[r]^{f^*} \ar[d] &
F\text{Perv} (\spec (A)^\N,\O_{\la\bullet}) \otimes L_\la \ar[d]
 \\
\mathbb{D}^b F (X^\N,\O_{\la\bullet})\otimes L_\la \ar[r]^{f^*} 
& 
\mathbb{D}^b F (\spec(A)^\N,\O_{\la\bullet})\otimes L_\la
}
$$ 
ce qui est immédiat.
\qed

\begin{lemm}
Dans le cas d'un point, il y a des équivalences de catégories 
$$
\mathbb{D}^b_c F (s,\O_\la)\otimes L_\la \iso \mathbb{D}^b_c F
(\O_\la) \otimes L_\la \iso \mathbb{D}^b_c F (L_\la)
$$
où $ \mathbb{D}^b_c F
(\O_\la)$ est la catégorie dérivée filtrée bornée des $\O_\la$-module
de type fini, $\mathbb{D}^b_c F (L_\la)$ celle des $L_\la$-e.v. de
dimension finie. La première équivalence est donnée par $R\pi_*$ et un
inverse par $\mathbb{L} \pi^*$.
\end{lemm}
\dem
Elle ne pose pas de problèmes particuliers
\qed

Désormais nous identifierons les trois catégories intervenant dans le
lemme précédent. 

\begin{prop}
Soit $\kappa : s \ldrt \spec (A)_s$ le point fermé. \`A isomorphisme
canonique près il existe un unique complexe $L(\F)\in \mathbb{D}^b_c F
(L_\la)$ tel que via le foncteur de normalisation
$\mathbb{D}^b_c F (s,\O_\la) \otimes L_\la \ldrt \mathbb{D}^b F
(s^\N,\O_{\la\bullet})\otimes L_\la$
$$
L(\F)\longmapsto \kappa^* K(\F)
$$
\end{prop}
\dem 
L'existence résulte de la proposition \ref{kitulo} et l'unicité de la
pleine fidélité du foncteur de normalisation.
\qed

Venons en au théorème principal.

\begin{theo}
Pour tout quadruplet $(X,x,\G,f)$ il y a un isomorphisme canonique induit par $f$ 
$$
K(\G)_x \iso L (\F)
$$
\end{theo}
\dem
Il s'agit de spécialiser la proposition \ref{kitulo} au point fermé de
$\spec (A)$.
\qed

La proposition suivante sera utilisée par Pascal Boyer dans son
travail. 

\begin{prop}
Soit $B=A [[ X_1,\dots, X_n]]$ et $h: \spec (B) \ldrt \spec (A)$ la
projection. Alors, 
$$
K (h^* \F [n] ) = h^* K ( \F) [n]
$$
$$
L ( \F) = L( h^* \F [n]) [-n]
$$
\end{prop}
\dem
Nous laissons les détails de la démonstration au lecteur. Il s'agit 
d'utiliser le théorème d'acyclicité local des
morphismes lisses couplé à la proposition \ref{kitulo} et à la
t-exactitude de $h^*$ (t-exactitude qui se déduit du théorème de
Popescu (dans ce cas là c'est un théorème de M.Artin) comme 
dans la démonstration du théorème \ref{cleft}). 
\qed

\bibliographystyle{smfplain}
\bibliography{biblio}

\providecommand{\bysame}{\leavevmode ---\ }
\providecommand{\og}{``}
\providecommand{\fg}{''}
\providecommand{\smfandname}{\&}
\providecommand{\smfedsname}{\'eds.}
\providecommand{\smfedname}{\'ed.}
\providecommand{\smfmastersthesisname}{M\'emoire}
\providecommand{\smfphdthesisname}{Th\`ese}
\begin{thebibliography}{10}

\bibitem{Berk3}
{\scshape V.~Berkovich} -- {\og Vanishing cycles for formal schemes. {II}\fg},
  \emph{Invent. Math.} \textbf{125} (1996), no.~2, p.~367--390.

\bibitem{BBD}
{\scshape J.~Bernstein, A.~Beilinson {\normalfont \smfandname} P.~Deligne} --
  {\og Faisceaux pervers\fg}, in \emph{Analyse et topologie sur les espaces
  singuliers, Asterisque 100}, 1982.

\bibitem{Boyer2}
{\scshape P.~Boyer} -- {\og Monodromie du faisceau pervers des cycles
  évanescents de quelques variétés de {S}himura simples et applications\fg},
  \emph{preprint}.

\bibitem{Eke}
{\scshape T.~Ekedahl} -- {\og On the adic formalism\fg}, in \emph{The
  {G}rothendieck {F}estschrieft}, vol.~2, Birkhauser.

\bibitem{FujGabber}
{\scshape K.~Fujiwara} -- {\og A proof of the absolute purity conjecture (after
  {G}abber)\fg}, in \emph{Algebraic geometry 2000, Azumino}, Adv. Stud. Pure
  Math., 36, 2002, p.~153--183.

\bibitem{Fuj2}
{\scshape K.~Fujiwara} -- {\og Theory of tubular neighborhood in \'etale
  topology\fg}, \emph{Duke Math. J.} \textbf{80} (1995), no.~1, p.~15--57.

\bibitem{Gabber1}
{\scshape O.~Gabber} -- {\og Notes on some {$t$}-structures\fg}, in
  \emph{Geometric aspects of Dwork theory. Vol. I, II}, Walter de Gruyter GmbH
  \& Co. KG, Berlin, 2004, p.~711--734.

\bibitem{Hu3}
{\scshape R.~Huber} -- {\og A finiteness result for direct image sheaves on the
  étale site of rigid analytic varieties\fg}, \emph{J. Algebraic Geom.}
  \textbf{7} (1998), no.~2, p.~159--403.

\bibitem{Illusie1}
{\scshape L.~Illusie} -- \emph{Complexe cotangent et déformations. {I}},
  Lecture Notes in Mathematics, vol. 239, Springer-Verlag, Berlin-New York,
  1971.

\bibitem{Illusie3}
\bysame , \emph{Complexe cotangent et déformations. {II}}, Lecture Notes in
  Mathematics, vol. 283, Springer-Verlag, Berlin-New York, 1972.

\bibitem{Jou1}
{\scshape J.-P. Jouanolou} -- {\og Cohomologie $\ell$-adique\fg}, in \emph{SGA
  V}.

\bibitem{Pop}
{\scshape M.~Spivakovsky} -- {\og A new proof of {D}. {P}opescu's theorem on
  smoothing of ring homomorphisms\fg}, \emph{J. Amer. Math. Soc.} \textbf{12}
  (1999), no.~2, p.~381--444.

\end{thebibliography}
\end{document}